\documentclass[12pt]{article}

\usepackage{amsmath}
\usepackage{amssymb}
\usepackage{amscd}
\usepackage{amsthm}
\theoremstyle{plain}
\newtheorem{Thm}{Theorem}[section]

\newtheorem{Prop}[Thm]{Proposition}

\newtheorem{Lemm}[Thm]{Lemma}
\newtheorem{Clai}[Thm]{Claim}
\theoremstyle{definition}
\newtheorem{Defn}[Thm]{Definition}
\theoremstyle{remark}
\newtheorem{Rem}[Thm]{Remark}
\numberwithin{equation}{section}
\def\Exc{\mathop{\mathrm{Exc}}}
\def\Rat{\mathop{\mathrm{Rat}}}

\def\Hilb{\mathop{\mathrm{Hilb}}}

\begin{document}
\title{Algebraic fiber space \\ whose generic fiber and base space are \\ of almost general type}
\author{ Shigetaka FUKUDA}
 
\date{\empty}

\maketitle \thispagestyle{empty}
\pagestyle{myheadings}
\markboth{Shigetaka FUKUDA}{Algebraic fiber space}
\begin{abstract}
We assume that the existence and termination conjecture for flips holds.
A complex projective manifold is said to be {\it of almost general type} if
the intersection number of the canonical divisor with every very general
curve is strictly positive.
Let $f$ be an algebraic fiber space from $X$ to $Y$.
Then the manifold $X$ is of almost general type if every very general fiber
$F$ and the base space $Y$ of $f$ are of almost general type.
\end{abstract}

\section{Introduction}

In this paper, every algebraic variety is defined over the field
$\mathbf{C}$ of complex numbers.
The two main conjectures of the classification theory of algebraic varieties
are the minimal model conjecture and the abundance conjecture.
The first consists of the existence and the termination of flips.
This conjecture appeared to be very natural, owing to the recent progress by
Shokurov \cite{Sh} and Hacon-McKernan \cite{HM}.
Furthermore the existence is now a theorem due to Birkar, Cascini, Hacon and
McKernan \cite{BCHM}.
The second is the conjecture that every nef canonical divisor should be
semi-ample.
But this is never trivial even for surfaces, as was seen in known proofs for
them.
Recently Ambro (\cite{Am}) reduced this abundance conjecture to the log
minimal model conjecture and to the conjecture
that every quasi-nup log canonical divisor should be semi-ample.
Here we explain concepts related to ``quasi-nup''.

%%%Def of almost nup and quasi-nup***************
\begin{Defn}
A $\mathbf{Q}$-Cartier divisor $D$ on a projective variety $X$ is {\it
almost numerically positive} ({\it almost nup}, for short), if there exists
a union $F$ of at most countably many prime divisors on $X$ such that $(D ,
C) > 0$ for every curve $C \nsubseteq F$ (i.e.\ if $(D , C) > 0$ for every
very general curve $C$).
We say that $D$ is {\it quasi-numerically positive} ({\it quasi-nup}, for
short) if $D$ is nef and almost nup.
\end{Defn}

%%%Def of almost general type*******************
\begin{Defn}
An algebraic variety $X$ is {\it of almost general type}, if there exists a
projective variety $M$ with only terminal singularities such that the
canonical divisor $K_M$ is almost nup and that $M$ is birationally equivalent to $X$.
\end{Defn}

\begin{Rem}
In the definition above, the assumption that $M$ is with only terminal singularities is important.
Indeed, there exists some ruled surface $X$ with wild singularities such that $K_X$ is big.
In this case, $X$ is not of almost general type but $K_X$ is almost nup.
See \cite{Ko} Example 43.
\end{Rem}

\begin{Prop}[Fukuda \cite{Fuk2}]\label{Prop:Fuk2 Douchi}
Under the minimal model conjecture, a projective variety with only terminal
singularities is of almost general type if and only if the canonical divisor
is almost nup.
\end{Prop}

\begin{Prop}[Fukuda \cite{Fuk2}]\label{Prop:MMCAC}
Under the minimal model conjecture and the abundance conjecture, a
projective variety is of almost general type if and only if it is of general type.
\end{Prop}

As was seen above, it is meaningful to research the class of varieties of
almost general type, depending not on the abundance conjecture but on the
minimal model conjecture (cf.\ \cite{Fuk} and \cite{Fuk2}).
Now we state the following

\begin{Thm}[Main Theorem]\label{Thm:MT}
Let $f \colon X \to Y$ be an algebraic fiber space such that every very general
fiber is of almost general type and that the base space $Y$ is of almost general type.
We assume that the minimal model conjecture, which consists of the existence
and the termination of flips, holds in dimension $\le \dim X$.
Then $X$ is of almost general type.
\end{Thm}

\begin{Rem}
By Proposition \ref{Prop:MMCAC}, if the abundance conjecture holds, then the theorem above is included in Viehweg (\cite{Vi}, Theorem 3), but the proof given in this paper is geometric.
\end{Rem}

\section{Preliminaries}

First we prepare basic notions and results concerning fibrations.

\begin{Defn}\label{Defn:fib}
A surjective morphism $f \colon X \to Y$ between projective varieties is {\it an
algebraic fiber space} if the extension $\Rat X / \Rat Y$ of a field is algebraically closed.
For a surjective morphism $f \colon X \to Y$ between projective varieties, we have
the naturally defined algebraic fiber space $[f]$ from the normalization of
$X$ to the normalization of $Y$ in $\Rat X$, by the Stein factorization.
\end{Defn}

\begin{Defn}
For an algebraic variety $X$, a subset $U$ is said to be {\it an almost open
dense subset} of $X$ if $U = X \setminus \bigcup_{i=1}^{\infty} V_i$ for
some countably many Zariski-closed subsets $V_i \subsetneqq X$.
The points that belong to some fixed almost open dense subset are said to be
{\it very general}.
\end{Defn}

\begin{Lemm}\label{Lemm:C}
If $f \colon X \to Y$ is a surjective morphism between complete algebraic
varieties and if $S$ is an almost open dense subset of $X$, then $f (S)$
contains some almost open dense subset of $Y$.
\end{Lemm}

%%%Proof of Lemma [C]***
\begin{proof}
First we treat the case where $f$ is an algebraic fiber space.
For some open dense subset $U'$ of $Y$, every fiber of $f \vert_{f^{-1}
(U')}$ is irreducible.
Take countably many open dense subsets $S_i$ of $X$ such that $S =
\bigcap_{i=1}^{\infty} S_i$.
From the theory of constructible sets, $f(S_i)$ contains some open dense
subset of $Y$.
Thus the set $U'' := \bigcap_{i=1}^{\infty} f(S_i) \cap U'$ contains some
almost open dense subset of $Y$.
For every point $y \in U''$, the intersection $f^{-1} (y) \cap S = f^{-1} 
(y) \cap \bigcap_{i=1}^{\infty} S_i$ is nonempty because $f^{-1} (y)$ is
irreducible.
Thus $f(S)$ contains $U''$ and therefore some almost open dense subset of $Y$.

Next we treat the case where $f$ is generically finite.
Then the image of every proper Zariski-closed subset of $X$ by $f$ is a
proper Zariski-closed subset of $Y$.
Thus the assertion follows in this case.

Lastly we treat the general case.
By considering the algebraic fiber space $[f]$ as in Definition \ref{Defn:fib}, the preceding two cases imply the assertion.
\end{proof}

\begin{Lemm}\label{Lemm:A}
Let $f \colon X \to Y$ be an algebraic fiber space and let $\{ V_t \}$ be an algebraic family of subvarieties of $X$.
Then one of the following two alternative possibilities occurs:

(i) $f(V_{t'})$ is a point for every general $t'$.

(ii) $f(V_{t'})$ is not a point for every general $t'$.
\end{Lemm}

%%%Proof of [A]***
\begin{proof}
By considering the intersection number $L^{m-1} \cdot f^* H \vert_{V_{t'}}$ for
some fixed ample divisors $L$ on $X$ and $H$ on $Y$ where $m= \dim V_{t'}$, we
obtain the assertion.
\end{proof}

The property ``of almost general type'' should be stable under taking a covering (Lemma \ref{Lemm:Bou}) and under deformations (Proposition \ref{Prop:Fuk2 Stab 
Deform}) and should have two geometric characterizations (Proposition \ref{Prop:Fuk2 MT}).

\begin{Lemm}\label{Lemm:Bou}
Let $f \colon X \to Y$ be a generically finite surjective morphism between
algebraic varieties.
If $Y$ is of almost general type, then also $X$ is of almost general type.
\end{Lemm}

%%%Proof of Lemma[Bou]***
\begin{proof}
By changing $X$ and $Y$ by adequate birational models, we may assume that
$X$ and $Y$ are nonsingular projective and that $K_{Y}$ is almost nup.
Thus the ramification formula implies the assertion.
\end{proof}

%%%[Aru Stab deform]***
\begin{Prop}[Fukuda \cite{Fuk2}]\label{Prop:Fuk2 Stab Deform}
Assume that the minimal model conjecture holds in dimension $n$.
Let $f \colon X \to Y$ be an algebraic fiber space with relative dimension $n$.
Then one of the following holds:

(i) Every general fiber $F$ of $f$ is not of almost general type.

(ii) Every very general fiber $F$ of $f$ is of almost general type.
\end{Prop}

\begin{Defn}
For an ambient space $X$, a subset $L$ is {\it covered} (resp.\ {\it swept out}) by subsets $D_i$ ($i
\in I$) if $L \subset \bigcup_{i \in I} D_i$ (resp.\ $L = \bigcup_{i \in I} D_i$).
\end{Defn}

%%%[Fuk2 Main Thm]*******************
\begin{Prop}[Fukuda \cite{Fuk2}]\label{Prop:Fuk2 MT}
Assume that the minimal model conjecture holds in dimension $\leq n$.
Let $X$ be a projective variety with only terminal singularities of dimension $n$.
Then the three conditions below are equivalent to each other:

(1) $X$ is of almost general type.

(2) The locus $\bigcup \{ D ; \thickspace D$ is a closed subvariety
($\subsetneqq X$) not of almost general type$\} $ is covered by at most
countably many prime divisors on $X$ and the variety $X$ is not birationally
equivalent to any minimal variety with numerically trivial canonical
divisor and is not a rational curve.

(3) $X$ is not uniruled and can not be covered by any family of varieties
being birationally equivalent to minimal varieties with numerically trivial
canonical divisors.
\end{Prop}

The following two lemmas are some kind of relativization for Proposition 
\ref{Prop:Fuk2 MT}.

\begin{Lemm}\label{Lemm:B}
Let $f \colon X \to Y$ be an algebraic fiber space with relative dimension $n$ 
such that every very general fiber is of almost general type.
Assume that the minimal model conjecture holds in dimension $\le n$.
Then any algebraic family of subvarieties that contract to points by $f$ and
that are not of almost general type can not cover $X$.
\end{Lemm}

%%%Proof of [B]***
\begin{proof}
Assume that such a family $\{ V_t \} $ covers $X$.
Then $\dim V_{t'} \leq [$the relative dimension of $f] = n$ for every
general member $V_{t'}$, because $f(V_{t'})$ is a point from Lemma \ref{Lemm:A} and the set $\{ f(V_{t'}) \}$ covers some open dense subset of $Y$.
Thus Proposition \ref{Prop:Fuk2 Stab Deform} implies that $V_{t'}$ is not of almost general type.
Therefore some open dense subset of every general fiber $F$ of $f$ is covered by some members $V_h$ ($\subset F$) not of almost general type.
Hence the fiber $F$ becomes not of almost general type from Proposition \ref{Prop:Fuk2 MT}.
This is a contradiction!
\end{proof}

\begin{Lemm}\label{Lemm:D}
Assume that the minimal model conjecture holds in dimension $\leq n$.
Let $f \colon X \to Y$ be an algebraic fiber space with relative dimension $n$.
Suppose that any algebraic family of subvarieties that contract
to points by $f$ and that are not of almost general type can not cover $X$.
Then the union of proper subvarieties that contract to points by $f$ and are
not of almost general type is contained in some union of at most countably
many proper Zariski-closed subsets of $X$.
\end{Lemm}

%%%Proof of Lemma[D]***
\begin{proof}
We follow the arguments in \cite{Fuk2} and, assuming that the locus $\bigcup \{ S ; \thickspace S$ is a closed subvariety ($\subsetneqq X$) not of almost general type and $f(S)$ is a point$\}$
cannot be covered by at most countably many prime divisors on $X$, derive a contradiction.

Let $\mathcal{H} \subset X \times \Hilb (X)$ be the universal family
parametrized by the Hilbert scheme $\Hilb (X)$.
By the countability of the components of $\Hilb (X)$, we have an irreducible
component $V$ of $\mathcal{H}$ with
surjective projection morphisms $\mu \colon V \to X$ and $\nu \colon V \to T ( \subset 
\Hilb (X))$ from $V$ to projective varieties
$X$ and $T$ respectively, such that $\mu(\nu ^{-1} (t)) \subsetneqq X$ for
every $t \in T$ and that the locus $\bigcup \{ S ; \thickspace S$ is a
closed subvariety ($\subsetneqq X$) not of almost
general type, $f(S)$ is a point and $S = \mu(\nu ^{-1} (t))$ for
some $t \in T \}$ cannot be covered by at most countably many prime divisors on $X$.
Let $\rho \colon V_{norm} \to V$ be the normalization.
We consider the Stein factorization of $\nu \rho$ into the finite morphism
$\nu _1 \colon W \to T$ from a projective normal
variety $W$ and the morphism $[\nu] \colon V_{norm} \to W$ with an algebraically
closed extension $\Rat V / \Rat W$.
Put $V^* :=$ [the image of the morphism $(\mu \rho, [\nu]) \colon V_{norm} \to X
\times W$].
\begin{equation}
\begin{CD}
V_{norm} @>>> V^* @>\text{embedding}>> X \times W @>>> W    \\
@.   @VVV @VVV @VV\text{$\nu_1$}V \\
@.   V   @>\text{embedding}>> X \times T @>>> T    \\
@.   @.   @VVV \\
@.   @.   X
\end{CD}
\end{equation}
Note that every fiber of the morphism $\nu \colon V \to T$ consists of a finite
number of fibers of the projection morphism from $V^*$ to $W$.
Thus we may replace $(V,T)$ by $(V^* , W)$ and assume that the extension
$\Rat V / \Rat T$ is algebraically closed.

According to Lemma \ref{Lemm:A} and Proposition \ref{Prop:Fuk2 Stab Deform} under the minimal model conjecture, we divide the situation into two cases.

First consider the case where $f(\mu(\nu^{-1} (t)))$ is a point (thus $\dim \mu (\nu^{-1} (t)) \le
[$the relative dimension of $f] = n$) and $\mu (\nu^{-1} (t))$ is of almost general type 
for very general $t \in T$ or where $f(\mu(\nu^{-1} (t)))$ is not a point for general $t \in T$.
Then there exists a subvariety $T_0 \subsetneqq T$ such that the locus
$\bigcup \{ S ; \thickspace S$ is a closed
subvariety ($\subsetneqq X$) not of almost general type, $f(S)$ is a point
and $S = \mu(\nu ^{-1} (t))$ for some $t \in T_0 \}$ cannot be covered by at
most countably many prime divisors on $X$.
Thus we can replace $(V, T)$ by $(V_1 , T_1)$, where $V_1$ and $T_1$ are
projective varieties such that $V_1$ is some suitable irreducible component
of $\nu ^{-1} (T_0)$ and $T_1 = \nu (V_1)$.
Because $\dim V_1 < \dim V$, by repeating this process of replacement, we
can reduce the assertion to the next case.

Now consider the case where $f(\mu(\nu^{-1} (t)))$ is a point and $\mu(\nu^{-1} (t))$ is not of almost general type for general $t \in T$.
This family $\{ \mu (\nu ^{-1} (t)) \}$ covers $X$.
So we obtain a contradiction!
\end{proof}
%%%*****************************************************

\section{Proof of Main Theorem}

Firstly we show that $X$ is not birationally equivalent to any minimal
projective variety $M$ with numerically trivial canonical divisor $K_M$.
Assuming that $X$ is birationally equivalent to such a variety $M$, we shall 
derive a contradiction.
We may suppose that $M$ is $\mathbf{Q}$-factorial, by executing $\mathbf{Q}$-factorialization.
Take a common resolution $X'$ of $X$ and $M$ with projective morphisms $g \colon X' 
\to X$ and $\mu \colon X' \to M$ such that $\Exc (\mu)$ is a divisor with only simple normal crossings.
Let $E$ be a $\mu$-exceptional prime divisor on $X'$ such that $f(g(E)) = Y$.
From Kawamata \cite{Kaw2}, the divisor $E$ is swept out by a family $\{ 
C_{\lambda} \}$ of $\mu$-exceptional rational curves.
The curves $C$ on $X'$ which are $fg$-exceptional are characterized by the
property that $((fg)^* H, C) = 0$, for some fixed ample divisor $H$ on $Y$.
Hence every general member $C_{\lambda'}$ of this family is a $fg$-exceptional 
rational curve, because $Y$ is not
uniruled by Miyaoka-Mori \cite{MiMo} from the assumption that $Y$ is of
almost general type.
Thus the covering family $\{ C_{\lambda} \}$ for $E$ is a family  of curves 
which are $\mu$-exceptional and $fg$-exceptional.
We consider a general fiber $F$ of $fg$.
Then $E \cap F$ consists of a finite number of irreducible fibers of the 
algebraic fiber space $[fg \vert_E]$.
Some nonempty open subset of every general (thus irreducible) fiber of $[fg 
\vert_E]$ is covered by general members $C_{\lambda'}$.
Therefore some open dense subset $U$ of $E \cap F$ is covered by sets $C_h \cap F$ ($\ne \emptyset$), where $C_h$ are some curves on $E$ which are $\mu$-exceptional and $fg$-exceptional.
Here $C_h$ is contained in $F$, because it is $fg$-exceptional, and therefore
$U$ is covered by these curves $C_h$ ($\subset E \cap F$).
Hence the divisor $E \vert_F$ is $\mu \vert_F$-exceptional.
Consequently the canonical divisor $K_F$ on $F$ is numerically equivalent to a $\mu \vert_F$-exceptional $\mathbf{Q}$-divisor, because $K_{X'}$ is numerically equivalent to a $\mu$-exceptional $\mathbf{Q}$-divisor.
So $K_F A^{\dim F -1} = 0$ for the pull-back $A$ of some ample divisor on $\mu (F)$.
Thus $K_F$ is not almost nup and therefore $F$ is not of almost general type 
by Proposition \ref{Prop:Fuk2 Douchi}.
This is a contradiction!

Secondly we prove that every very general member of any non-trivial 
algebraic family that covers $X$ is of almost general type.
By Proposition \ref{Prop:Fuk2 Stab Deform}, it suffices to derive a 
contradiction, assuming that every general member of some non-trivial 
algebraic family $\{ S_t \}_{t \in T}$ that covers $X$ is not of almost general type.

\begin{Clai}\label{Clai:Claim}
For some member $D$ of $\{ S_t \}_{t \in T}$, we have that $f(D)=E$ and $1 \leq \dim E \leq \dim D < \dim X$, where $D$ is not of almost general type and $E$ is of almost general type and where every very general fiber of the algebraic fiber space $[f \vert_D ]$ is of almost general type or the morphism $f \vert_D$ is generically finite.
\end{Clai}

%%%Proof of Claim******
\begin{proof}
We note that $S_t \subsetneqq X$ for general $t$, because $\{ S_t \}$ is a non-trivial algebraic family.

From Lemmas \ref{Lemm:A} and \ref{Lemm:B}, the image $f(D)$ is not a point
for every general member $D$ of $\{ S_t \}$.

For every almost open dense subset $U$ of $Y$, the set $\{ t \in T ; f^{-1} (U) \cap S_t \ne \emptyset \}$ contains some almost open dense subset of $T$, from Lemma \ref{Lemm:C}.
So the variety $f(D)$ becomes very general subvariety of $Y$ for every very general
member $D$ of $\{ S_t \}$ and therefore it is of almost general type from Proposition 
\ref{Prop:Fuk2 MT}.

Any algebraic family of subvarieties that contract to points by
$f$ (cf.\ Lemma \ref{Lemm:A}) and that are not of
almost general type (cf.\ Proposition \ref{Prop:Fuk2 Stab Deform}) can
not cover $X$, from Lemma \ref{Lemm:B}.
Thus the union of proper subvarieties that contract to points by $f$ and are
not of almost general type is contained in some union of at most countably
many proper Zariski-closed subsets of $X$, from Lemma \ref{Lemm:D}.
Consequently every very general fiber of $[f \vert_D]$ is of almost general
type or the morphism $f \vert_D$ is generically finite.

The proof for Claim \ref{Clai:Claim} ends.
\end{proof}

\noindent
In the situation of Claim \ref{Clai:Claim}, by applying the induction hypothesis for $(X,Y)$ on $\dim X$ to $(D,E)$ and by Lemma \ref{Lemm:Bou}, we obtain the conclusion that $D$ is of almost general type.
But this contradicts the statement that $D$ is not of almost general type!

By Proposition \ref{Prop:Fuk2 Douchi}, every minimal variety of dimension $\le \dim X$ with 
numerically trivial canonical divisor is not of almost general type.
Thus, from the first and second paragraph of this section, $X$ is not covered by any 
algebraic family of varieties being
birationally equivalent to minimal varieties with numerically trivial canonical divisors.
Because the rational curve is not of almost general type, $X$ is also not 
uniruled from the second paragraph of this section.
Consequently Proposition \ref{Prop:Fuk2 MT} implies the assertion.
\qed

\noindent{\bf Acknowledgement}

\noindent
The author would like to thank the referees who read the manuscript carefully, gave many valuable advice to improve the presentation, and informed him of the relevant example.

\bigskip
Faculty of Education, Gifu Shotoku Gakuen University

Yanaizu-cho, Gifu City, Gifu 501-6194, Japan

fukuda@ha.shotoku.ac.jp

\end{document}